\documentclass{article}[10pt]
\usepackage[utf8]{inputenc}
\usepackage{
 amsmath,
 amsxtra,
 amsthm,
 amssymb,
 etex,
 mathrsfs,
  }
\usepackage[all]{xy}
\usepackage[T2A,T1]{fontenc}
\usepackage[OT2,T1]{fontenc}
\newcommand{\Zhe}{\mbox{\usefont{T2A}{\rmdefault}{m}{n}\CYRZH}}

\newtheorem{theorem}{Theorem}[section]
\newtheorem{lemma}[theorem]{Lemma}

\newtheorem{proposition}[theorem]{Proposition}
\newtheorem{corollary}[theorem]{Corollary}

\theoremstyle{definition}
\newtheorem{defn}[theorem]{Definition}
\newtheorem{remark}[theorem]{Remark}

\newcommand{\bd}{\begin{defn}}
\newcommand{\ed}{\end{defn}}
\newcommand{\bl}{\begin{lemma}}
\newcommand{\el}{\end{lemma}}
\newcommand{\bp}{\begin{proposition}}
\newcommand{\ep}{\end{proposition}}
\newcommand{\bt}{\begin{theorem}}
\newcommand{\et}{\end{theorem}}
\newcommand{\bc}{\begin{corollary}}
\newcommand{\ec}{\end{corollary}}
\newcommand{\br}{\begin{remark}}
\newcommand{\er}{\end{remark}}
\newcommand{\ba}{\begin{array}}
\newcommand{\ea}{\end{array}}
\newcommand{\bpf}{\begin{proof}}
\newcommand{\epf}{\end{proof}}

\newcommand{\Z}{\mathbb{Z}}
\newcommand{\Q}{\mathbb{Q}}
\newcommand{\Zp}{\mathbb{Z}_{p}}
\newcommand{\Qp}{\mathbb{Q}_{p}}

\newcommand{\Ap}{A[p^{\infty}]}

\newcommand{\Ga}{\Gamma}

\newcommand{\La}{\Lambda}
\newcommand{\la}{\lambda}

\DeclareMathOperator{\Sel}{Sel} \DeclareMathOperator{\Gal}{Gal}
 \DeclareMathOperator{\rank}{rank}

 \DeclareMathOperator{\Tor}{Tor}

\newcommand{\cyc}{\mathrm{cyc}}
\newcommand{\ac}{\mathrm{ac}}

\newcommand{\G}{\mathfrak{G}}

\newcommand{\mM}{\mathcal{M}}

\newcommand{\ot}{\otimes}
\newcommand{\ilim}{\displaystyle \mathop{\varinjlim}\limits}
\newcommand{\plim}{\displaystyle \mathop{\varprojlim}\limits}

\newcommand{\coker}{\mathrm{coker}\,}

\newcommand{\lra}{\longrightarrow}
\newcommand{\tha}{\twoheadrightarrow}

\newcommand{\ps}[1]{[[ #1 ]]}

  \DeclareFontFamily{U}{wncy}{}
  \DeclareFontShape{U}{wncy}{m}{n}{<->wncyr10}{}
  \DeclareSymbolFont{mcy}{U}{wncy}{m}{n}
  \DeclareMathSymbol{\sha}{\mathord}{mcy}{"58}

\numberwithin{equation}{subsection}

\begin{document}
\title{On the control theorem for fine Selmer groups and the growth of fine Tate-Shafarevich groups in $\Zp$-extensions}
 \author{Meng Fai Lim\footnote{School of Mathematics and Statistics $\&$ Hubei Key Laboratory of Mathematical Sciences, Central China Normal University, Wuhan, 430079, P.R.China. E-mail: \texttt{limmf@mail.ccnu.edu.cn}} }
\date{}
\maketitle

\begin{abstract} \footnotesize
\noindent Let $A$ be an abelian variety defined over a number field $F$. We prove a control theorem for the fine Selmer group of the abelian variety $A$ which essentially says that the kernel and cokernel of the natural restriction maps in an arbitrarily given $\Zp$-extension $F_\infty/F$ are finite and bounded. We emphasise that our result does not have any constraints on the reduction of $A$ and the ramification of $F_\infty/F$. As a first consequence of the control theorem, we show that the fine Tate-Shafarevich group over an arbitrary $\Zp$-extension has trivial $\La$-corank. We then derive an asymptotic growth formula for the $p$-torsion subgroup of the dual fine Selmer group in a $\Zp$-extension. However, as the fine Mordell-Weil group needs not be $p$-divisible in general, the fine Tate-Shafarevich group needs not agree with the $p$-torsion of the dual fine Selmer group, and so the  asymptotic growth formula for the dual fine Selmer groups do not carry over to the fine Tate-Shafarevich groups. Nevertheless, we do provide certain sufficient conditions, where one can obtain a precise asymptotic formula.


\end{abstract}

\section{Introduction}

The essence of Iwasawa theory is to study arithmetic objects via their variations in a tower of number fields. Such a study was initiated by Iwasawa, where he
showed that for a given prime $p$, the growth of the Sylow $p$-subgroups of the class groups in the tower of subfields of a $\Zp$-extension exhibits a remarkable regularity (see \cite{Iw59, Iw73}). Modelling after Iwasawa's ideas, Mazur developed an analogous approach towards studying the arithmetic of an abelian variety via examining the variations of its $p$-primary Selmer groups in a $\Zp$-extension (see \cite{Maz}). Recently, there have been great interest in the study of a certain subgroup of the $p$-primary Selmer group, called the fine Selmer group (for instances, see \cite{CS05,JS,LimFine,LMu,LSu,LSu2,Su10,Wu}). In the fundamental paper of Coates-Sujatha \cite{CS05}, they have examined this fine Selmer group in great depth and made several conjectures on its structure. We should also mention that (as also noted by Coates-Sujatha in \cite{CS05}) before the foundational work of Coates-Sujatha, this group has been studied in various guises (for instances, see \cite{K, Kob, Kur, PR00}).

Just as the classical $p$-primary Selmer group $\Sel_p(A/F)$ sitting in the middle of a short exact sequence
 \[0\lra A(F)\ot_{\Zp}\Qp/\Zp \lra \Sel_p(A/F) \lra \sha(A/F)[p^\infty] \lra 0,\]
with the Mordell-Weil group and the $p$-primary Tate-Shafarevich group by its sides, the fine Selmer group $R_p(A/F)$ (see Section \ref{fine selmer sec} for definition) sits in the middle of the following analogous short exact sequence
 \[0\lra \mM_p(A/F) \lra R_p(A/F) \lra \Zhe_p(A/F) \lra 0,\]
where $\mM_p(A/F)$ is the fine ($p$-)Mordell-Weil group and  $\Zhe_p(A/F)$  is the fine ($p$-)Tate-Shafarevich group in the sense of Wuthrich \cite{WuTS}. Here $\mM_p(A/F)$ is defined to be the subgroup of $A(F)\ot\Qp/\Zp$ consisting of those elements which are mapped to zero in $A(F_v)\ot_{\Zp}\Qp/\Zp$ for all primes $v$ above $p$. It is not difficult to show that $\mM_p(A/F)$ injects into $R_p(A/F)$. The fine Tate-Shafarevich group $\Zhe_p(A/F)$ is then simply defined to be the cokernel of this injection, and one can verify that $\Zhe_p(A/F)$ can be identified as a subgroup of $\sha(A/F)[p^\infty]$ (see \cite[Section 2]{WuTS}; also see discussion below in Section \ref{fine selmer sec}).  In particular, we would expect that the fine Tate-Shafarevich group is finite in view of the conjectural finiteness of $\sha(A/F)$. We also make a remark that will be frequently mentioned throughout the paper, namely that, the fine Mordell-Weil group $\mM_p(A/F)$ is not necessarily $p$-divisible (see \cite[Section 7]{WuTS} for examples of non-divisible fine Mordell-Weil groups).

\textbf{As our prime $p$ is fixed throughout, we shall omit the $``p"$ in our notations of the above arithmetic objects for the remainder of the paper. We shall also always assume that the prime $p$ is odd.} (See the end of the introductory section for some remarks on the case $p=2$.)

In this paper, we are interested in studying the variation of the fine Selmer groups and the fine Tate-Shafarevich groups of an abelian variety $A$ in a $\Zp$-extension $F_\infty$. We write $\Ga=\Gal(F_\infty/F)$ and $\Zp\ps{\Ga}$ for the Iwasawa algebra of $\Ga$. For each $n\geq 0$, $F_n$ will denote the intermediate subfield of $F_\infty$ with $|F_n:F| = p^n$. We write $R(A/F_\infty)$ for the fine Selmer group of $A$ over $F_\infty$, whose Pontryagin dual is denoted by $Y(A/F_\infty)$. It is not difficult to show that this is a finitely generated $\Zp\ps{\Ga}$-module (see Lemma \ref{fg La modules}). However, it is unknown at present
whether or not $Y(A/F_\infty)$ is always $\Zp\ps{\Ga}$-torsion, although we believe that this is the case (see Conjecture Y in the body of the paper).

 We now state the following control theorem for fine Selmer groups.

\bt[Theorem \ref{control theorem}]
Let $A$ be an abelian variety defined over a number field $F$. Let $F_{\infty}$ be a $\Zp$-extension of $F$ and $F_n$ the intermediate subfield of $F_\infty$ with $|F_n:F| = p^n$.
Then the restriction map
\[r_n: R(A/F_n) \lra R(A/F_{\infty})^{\Gal(F_\infty/F_n)}\]
has finite kernel and cokernel which are bounded independent of $n$.
\et

For the classical $p$-primary Selmer groups, such a control theorem was established by Mazur \cite{Maz} under the assumption that the abelian variety has good ordinary reduction at all primes of $F$ above $p$. For an elliptic curve defined with split multiplicative reduction at the prime $p$, one generally expects a control theorem for the classical Selmer group but so far this is only known when the said elliptic curve is defined over $\Q$ (see \cite{G99}). However, if the abelian variety in question has good supersingular reduction at some prime above $p$, it is well-known that such a control theorem fails for the classical Selmer groups (see \cite{Kob, Sch85}). Our theorem here is saying that this is not an issue for the fine Selmer groups. We should also mention that when $A$ has good reduction at every prime of $F$ above $p$ and $F_\infty$ is the cyclotomic $\Zp$-extension, a control theorem for the fine Selmer group has been established (for instances, see \cite{Kob, LLSha, LP, WuPhD, Wu}). Our result can therefore be seen as a generalization of these prior results, as we do not have any constraints on the reduction of $A$ and the $\Zp$-extension $F_\infty/F$. As a corollary, we show that $Y(A/F_\infty)$ is $\Zp\ps{\Ga}$-torsion whenever $Y(A/F)$ is finite (see Corollary \ref{control theorem cor}).

As one of the further applications of our control theorem, we establish the following interesting phenomenon (see Proposition \ref{torsion fine TateSha} for a more precise statement).

\bp[Proposition \ref{torsion fine TateSha}]
 Let $A$ be an abelian variety defined over a number field $F$ and $F_{\infty}$ a $\Zp$-extension of $F$ with intermediate subfield $F_n$. Suppose further that $\Zhe(A/F_n)$ is finite for every $n$. Then $\Zhe(A/F_{\infty})$ is a cotorsion $\Zp\ps{\Ga}$-module.
\ep

Note that the corresponding assertion for the $p$-primary Tate-Shafarevich group is known to be false in general. Indeed, for an elliptic curve defined over $\Q$ with supersingular reduction at the prime $p$, it has been long observed that its $p$-primary Tate-Shafarevich group over the cyclotomic $\Zp$-extension is a non-cotorsion $\Zp\ps{\Ga}$-module by a combination of results of Kato-Rohrlich \cite{K, Ro} and Schneider \cite{Sch85} (also see \cite{Kob, Kur}).

We now come to the second theme of the paper concerning about the growth of the fine Tate-Shafarevich groups in a $\Zp$-extension. In preparation of subsequent discussion, we introduce certain terminology. For a finitely generated
$\Zp$-module $N$, write $e(N)$ for the power of $p$ in the order of $N[p^\infty]$, i.e., $\big|N[p^\infty]\big| = p^{e(N)}$. Now, if $M$ is a finitely generated $\Zp\ps{\Ga}$-module $M$, we follow Lee \cite{Lee} in defining
\[\G(M): = \plim_n \Big(M_{\Ga_n}[p^\infty]\Big). \]
Note that $\G(M)$ is always a finitely generated torsion $\Zp\ps{\Ga}$-module and has the same $\mu$-invariant as $M$ (see Proposition \ref{G functor}).
The Pontryagin dual of the fine Selmer groups $R(A/F_n)$ and $R(A/F_\infty)$ are denoted by $Y(A/F_n)$ and $Y(A/F_\infty)$ respectively.

 We can now state the following which is an immediate consequence of our control theorem.

\bp[Proposition \ref{growthfineselmer}]
 Let $A$ be an abelian variety defined over a number field $F$. Let $F_{\infty}$ be a $\Zp$-extension of $F$ and $F_n$ the intermediate subfield of $F_\infty$ with $|F_n:F| = p^n$. Then we have
 \[ e\big(Y(A/F_n)\big) = \mu\Big(\G\big(Y(A/F_\infty)\big)\Big)p^n +\la\Big(\G\big(Y(A/F_\infty)\big)\Big)n + O(1).\]
\ep

Note that we \textit{do not} assume the torsionness of the fine Selmer group in the above proposition. In view of the above formula, it is natural to ask if one can transfer the above growth formula to the fine Tate-Shafarevich groups as done for the usual Tate-Shafarevich groups in \cite{G99, Lee}. In these said works, a crucial ingredient used is the $p$-divisibility of $A(F_n)\ot_{\Zp}\Qp/\Zp$ which is a consequence from its definition. However, as already mentioned above, the fine Mordell-Weil group \textit{needs not} be $p$-divisible, and so the groups $\Zhe(A/F_n)$ and $Y(A/F_n)[p^\infty]$ may not agree in general. Therefore, the above asymptotic formula for the fine Selmer groups does not carry over to the fine Tate-Shafarevich groups directly.

We do at the very least have the following general asymptotic lower bound and upper bound.

\bp[Proposition \ref{growth bound}]
  Let $A$ be an abelian variety defined over a number field $F$ and $F_{\infty}$ a $\Zp$-extension of $F$. Suppose further that $\Zhe(A/F_n)$ is finite for every intermediate subfield $F_n$.
    Then we have
 \[ \mu\Big(\Zhe(E/F_\infty)^\vee\Big)p^n+  \la\Big(\Zhe(E/F_\infty)^\vee\Big)n + O(1) \leq e\big(\Zhe(A/F_n)\big) \hspace{2in}\]\[\hspace{2in} \leq \mu\Big(\G\big(Y(A/F_\infty)\big)\Big)p^n +\la\Big(\G\big(Y(A/F_\infty)\big)\Big)n + O(1).\]
\ep

Naturally, one might ask whether one can have an exact asymptotic growth formula under appropriate extra assumptions. The next theorem is the first of which where we can establish an exact asymptotic formula under the additional hypothesis that the Mordell-Weil group $A(F_\infty)$ over the $\Zp$-extension $F_\infty$ is finitely generated.

\bt[Theorems \ref{growth}]
 Let $A$ be an abelian variety defined over a number field $F$ and $F_{\infty}$ a $\Zp$-extension of $F$. Suppose that $A(F_\infty)$ is a finitely generated abelian group and that $\Zhe(A/F_n)$ is finite for every $n$.
    Then we have
 \[e\big(\Zhe(A/F_n)\big) = \mu\Big(\Zhe(A/F_\infty)^\vee\Big)p^n + \la\Big(\Zhe(A/F_\infty)^\vee\Big)n + \nu\]
 for $n\gg0$.
\et

For the cyclotomic $\Zp$-extension $F^\cyc$ of $F$, it has long been conjectured that $A(F^\cyc)$ is a finitely generated abelian group (see \cite{G99, Maz}). Til today, this finite generation conjecture has only been verified  by the deep works of Kato \cite{K} and Rohrlich \cite{Ro} under the assumptions that $A$ is an elliptic curve defined over $\Q$ and that $F$ is an abelian extension of $\Q$.  Therefore, we do have a precise asymptotic formula growth in this situation. This in particular recovers an assertion of Wuthrich \cite[Discussion before Conjecture 8.2]{WuTS} (also see \cite{WuPhD}).

On the other hand, over an arbitrary $\Zp$-extension, the group $A(F_\infty)$ needs not be finitely generated. Indeed, if $A$ is an elliptic curve over $\Q$ and $F_\infty$ is the anticyclotomic $\Zp$-extension of an imaginary quadratic field $F$, then the group $A(F_\infty)$ is expected to be not finitely generated under a certain root number condition (see \cite[Growth Number Conjecture]{Maz83} and \cite[Conjecture 1.2]{LSp}). This therefore brings us to the next theorem which provides us some criterion to have a precise growth formula in the absence of the finite generation of $A(F_\infty)$, but under more restrictive assumptions on the reduction type of $A$ and the finiteness of the full ($p$-primary) Tate-Shafarevich groups.

\bt[Theorems \ref{growth2}]
 Let $A$ be an abelian variety defined over a number field $F$ with potentially good ordinary reduction at all primes of $F$ above $p$. Let $F_{\infty}$ be a $\Zp$-extension of $F$ and suppose that $\sha(A/F_n)[p^\infty]$ is finite for every $n$. Then we have
 \[e\big(\Zhe(A/F_n)\big) = \mu\Big(\Zhe(A/F_\infty)^\vee\Big)p^n + \la\Big(\Zhe(A/F_\infty)^\vee\Big)n + \nu\]
 for $n\gg0$.
\et

We say a little on the proofs of Theorems \ref{growth} and \ref{growth2}. The ingredients of the proofs are our control theorem (Theorem \ref{control theorem}), a description of the structure of $\Zhe(A/F_\infty)^{\vee}$ via the $\G$-functor of Lee (see Proposition \ref{torsion fine TateSha}) and an explicit description of the structure of $(A(F_\infty)\ot_{\Zp}\Qp/\Zp)^\vee$ as elucidated in the work of Lee \cite{Lee}. Although at first viewing, Theorem \ref{growth2} does not have any prior hypothesis on the Mordell-Weil groups, we do made use of (a consequence of) Lee's explicit description for our eventual proof. We also mention that our approach in proving Theorem \ref{growth2} requires us to make use of the control theorem of the classical Selmer group of the ordinary abelian variety $A$ (but we do not require this ordinarity hypothesis for the proof of Theorem \ref{growth}).

We finally say something on the situation $p=2$. If the number field $F$ has no real primes, then all the above mentioned results carry over. In the event that $F$ has at least one real prime, the situation is slightly trickier due to technical cohomological considerations. Despite this, we are hopeful that some variant of the results in this paper should hold for $p=2$, although at this point of writing, we are not able to pinpoint the exact variant as yet. However, we like to mention that there are some recent interesting works \cite{CKL, Ke} on the Iwasawa theory for (possibly non-cyclotomic) $\Z_2$-extensions. Naturally, one might ask if the ideas in these works can be applied to study control theorem of fine Selmer groups and growth of fine Tate-Shafarevich groups for the case $p=2$. We hope to revisit this in the near future.

We now give an outline of the paper. In Section \ref{algebra prelim}, we collect several results on $\Zp$-modules and $\La$-modules which will be required in our arithmetic discussion. In Section \ref{fine selmer sec}, we introduce the fine Selmer groups, fine Mordell-Weil groups and fine Tate-Shafarevich groups. The control theorem of the fine Selmer groups will be established here. Section \ref{fine Zhe sec} is where we prove our asymptotic estimates on the growth of the fine Tate-Shafarevich groups in a $\Zp$-extension. In Section \ref{further comments}, we specialize our theorem to certain specific $\Zp$-extensions. In particular, for the cyclotomic $\Zp$-extension, we formalize a conjectural explicit growth of the fine Tate-Shafarevich groups which was first proposed by Wuthrich \cite{WuPhD, WuTS}. We also discuss some classes of non-cyclotomic $\Zp$-extensions, where we can obtain a precise growth formula of the fine Tate-Shafarevich groups. Along the way, we propose several questions on some finer aspects of the growth of fine Tate-Shafarevich groups in various $\Zp$-extensions.

\subsection*{Acknowledgement}
 The author would like to thank John Coates, Ming-Lun Hsieh, Antonio Lei and Christian Wuthrich for their interest and comments. He especially thanks Christian Wuthrich for discussion and answering questions on his works \cite{WuPhD, Wu, WuTS}. He also likes to thank Ming-Lun Hsieh for answering his questions pertaining to the anticyclotomic $\Zp$-extensions. The author is also grateful to the anonymous referee for many useful comments and suggestions on the article. Finally, he likes to acknowledge support by the National Natural Science Foundation of China under Grant No. 11550110172 and Grant No. 11771164.

\section{Algebraic preliminaries} \label{algebra prelim}

In this section, we recall certain algebraic preliminaries. If $N$ is a $\Zp$-module, denote by $N[p^n]$ the submodule of $N$ consisting of elements
of $N$ which are annihilated by $p^n$. We then write $N[p^\infty] = \cup_{n\geq 1}N[p^n]$. If $N$ is finitely generated
over $\Zp$, we write $e(N)$ for the power of $p$ in the order of $N[p^\infty]$, i.e., $|N[p^\infty]| = p^{e(N)}$.

\subsection{Some useful lemmas}
We begin with a useful lemma (compare with \cite[Lemma 2.1.4(2)]{Lee} and \cite[Lemma 2.1.3]{LiangL}).

\bl \label{short exact fg Zp-modules}
 Let $0\lra M \lra N \lra P\lra 0$ be a short exact sequence of finitely generated $\Zp$-modules. Suppose that $p^j$ annihilates $P[p^\infty]$. Then we have an exact sequence
 \[  0 \lra M[p^\infty] \lra N[p^\infty] \lra P[p^\infty] \lra M_f/p^j \]
 of abelian groups, where $M_f: = M/M[p^\infty]$.
 In particular, if $M$ is finite, we have a short exact sequence
 \[  0 \lra M[p^\infty] \lra N[p^\infty] \lra P[p^\infty] \lra 0 \]
 and an equality
 \[ e(N)= e(M)+ e(P).\]
 \el

\bpf
By considering the long exact $\Tor$-sequence of $-\ot_{\Zp}{\Qp/\Zp}$, we obtain an exact sequence
\[ 0 \lra  M[p^\infty] \lra N[p^\infty] \lra P[p^\infty] \stackrel{\delta}{\lra} M\ot\Qp/\Zp. \]
Since $p^j$ annihilates $P[p^\infty]$, the image of the map $\delta$ is contained in
\[ \big(M\ot\Qp/\Zp\big)[p^j] = \big(M_f\ot\Qp/\Zp\big)[p^j] =M_f/p^j.\]
This yields the exact sequence of the lemma. The remainder of the lemma is  immediate from this.
\epf

Occasionally, we need to analyze the $p$-torsion subgroups of the terms in an exact sequence with at least four terms. The following lemma gives us some leverage towards this.

\bl \label{exact fg Zp-modules}
 Let $0\lra M \lra N \lra P\lra Q $ be an exact sequence of finitely generated $\Zp$-modules, where $M$ is finite. Then we have an exact sequence
 \[  0 \lra M \lra N[p^\infty] \lra P[p^\infty] \lra Q[p^\infty] \]
 of $\Zp$-modules. In particular, we have
  \[ e(P)\leq e(N)+ e(Q) \quad \mbox{and}\quad \Big|e(N) - e(P) \Big| \leq e(M) + e(Q).\]
 \el

\bpf
 Write $U$ for the image of the map $N\lra P$. Then by Lemma \ref{short exact fg Zp-modules} and the finiteness of $M$, we have exact sequences
 \[
 0 \lra M \lra N[p^\infty]\lra U[p^\infty] \lra 0,
  \]\[ 0 \lra U[p^\infty] \lra P[p^\infty] \lra Q[p^\infty].
 \]
  Splicing the two exact sequences, we obtain the required exact sequence of the lemma. The estimates in the lemma are immediate from this.
\epf

As we frequently need to consider inverse limits of modules, the next lemma will come in handy for our subsequent discussion.

\bl \label{limit}
Let $\{M_n\}$, $\{N_n\}$, $\{P_n\}$ and $\{Q_n\}$ be projective systems of finitely generated $\Zp$-modules such that for every $n$, there is an exact sequence
 \[0\lra M_n \lra N_n \lra P_n\lra Q_n \]
 which is compatible with the transition maps of the projective systems. Suppose further that each $M_n$ is finite.
Then we have an exact sequence
 \[  0 \lra \plim_n M_n \lra \plim_n N_n \lra \plim_n P_n \lra \plim_n Q_n \]
 of $\Zp$-modules.
\el

\bpf
 Write $U_n$ for the image of the map $N_n\lra P_n$. Then we have exact sequences
\[
 0 \lra \plim_n M_n \lra \plim_n N_n \lra \plim_n U_n \lra 0,
  \]\[ 0 \lra \plim_n U_n \lra \plim_n P_n \lra \plim_n Q_n,
 \]
 where the surjectivity of the rightmost map of the first sequence follows from the fact that $\plim_n\!^1 M_n =0$ by the finiteness of $M_n$. The required exact sequence of the lemma then follows from splicing the above two exact sequences. \epf

\subsection{The functor $\mathfrak{G}$}

Throughout the paper, we shall write $\La$ for the classical Iwasawa algebra $\Zp\ps{\Ga}$, where $\Ga=\Zp$.  For a finitely generated $\La$-module $M$, we can attach Iwasawa $\mu$-invariant (denoted by $\mu(M))$ and Iwasawa $\la$-invariant (denoted by $\la(M))$ to it (see \cite[Definition 5.3.9]{NSW}).

Denote by $\Ga_n$ the unique subgroup of $\Ga$ of index $p^n$. In \cite{Lee}, Lee introduced the following functor on a $\La$-module $M$ which is defined by
\[\G(M): = \plim_n \Big(M_{\Ga_n}[p^\infty]\Big). \]
We shall record certain properties of this functor that will be required in our discussion. For more details on the functor, we refer the readers to the paper \cite{Lee}.

\bp \label{G functor}
Let $M$ be a finitely generated $\La$-module. Then the following assertions are valid.
\begin{enumerate}
\item[$(1)$] $\G(M)$ is a torsion $\La$-module.
\item[$(2)$] $\mu(\G(M)) = \mu(M)$.
\item[$(3)$] $\G(M)=0$ if and only if there is a pseudo-isomorphism
  \[M\stackrel{\sim}{\lra} \La^r \oplus \Big(\bigoplus_{j=1}^s \La/\xi_{m_j}\Big),\]
  where $\xi_{m_j}$ is certain $p^{m_j}$-th cyclotomic polynomial.
\item[$(4)$] We have \[e\big(M_{\Ga_n}\big) = \mu\big(\G(M)\big)p^n + \la\big(\G(M)\big)n + \nu\]
for $n\gg 0$.
\end{enumerate}
\ep

\bpf
Assertions (1)-(3) are immediate consequences of \cite[Lemma A.2.9]{Lee}.
 The final assertion is \cite[Lemma 4.1.3]{Lee}.
\epf

We end with the following technical result which will be required for the estimation of the growth of the fine Tate-Shafarevich groups.

\bl \label{tech estimate}
Let $\{M_n\}$ be a projective system of $\La$-modules with transition maps $M_{n+1}\lra M_n$ such that the action of $\La$ on $M_n$ factors through $\Zp[\Ga/\Ga_n]$. Write $M = \plim_n M_n$. For each $n$, the natural map $M\lra M_n$ factors through $M_{\Ga_n}$ to induce a map $M_{\Ga_n}\lra M_n$.
Suppose that the following statements are valid.
\begin{enumerate}
 \item[$(a)$] The kernel of the map $M_{\Ga_n}\lra M_n$ is finite for each $n$.
 \item[$(b)$] The cokernel of the map $M_{\Ga_n}\lra M_n$ is finite and bounded independently of $n$.
     \end{enumerate}
 Then we have $e(M_n) \leq \mu\big(\G(M)\big)p^n+\la\big(\G(M)\big)n+O(1)$. Furthermore, we have the following assertions.
 \begin{enumerate}
 \item[$(i)$] In the event that the kernel of the map $M_{\Ga_n}\lra M_n$ is also bounded independently of $n$, we then have $e(M_n) = \mu\big(\G(M)\big)p^n+\la\big(\G(M)\big)n+O(1)$.
 \item[$(ii)$] If $\G(M)= 0$, then the quantity $e(M_n)$ is bounded independently of $n$.
     \end{enumerate}
\el

\bpf
Denote by $C_n$ (resp., $D_n$) the kernel (resp., cokernel) of the map $M_{\Ga_n}\lra M_n$. Since $C_n$ is assumed to be finite, it follows from Lemma \ref{exact fg Zp-modules} that we have
\[ e(M_n) \leq e(M_{\Ga_n}) + e(D_n).\]
The required estimate then follows from Proposition \ref{G functor}(4) and hypothesis (b) of the proposition.
Now, if $\G(M)= 0$, then we see immediately that $e(M_n) = O(1)$. This establishes (ii). On the other hand, by Lemma \ref{exact fg Zp-modules} again, we have
\[ \big|e(M_n)- e(M_{\Ga_n})\big| \leq  e(C_n) + e(D_n).\]
Therefore, if the order of $C_n$ is also bounded, it follows from the above that we have $e(M_n) = \mu\big(\G(M)\big)p^n+\la\big(\G(M)\big)n+O(1)$. This shows (i).
\epf

\br
Note that in assertion (ii) of Lemma 2.5, we do not a priori require the hypothesis that the kernel of the map $M_{\Ga_n}\lra M_n$ is bounded independently of $n$.
\er

\section{Fine Selmer groups over $\Zp$-extension} \label{fine selmer sec}

\subsection{General setup}
We begin with some general remarks. Throughout, we let $A$ denote an abelian variety which is defined over a number field $F$. Let $S$ be a finite set of primes of $F$ containing the primes above $p$, the bad reduction primes of $A$ and the infinite primes. We then write $S_p$ for the set of primes in $S$ lying above $p$. Denote by $F_S$ the maximal algebraic extension of $F$ which is unramified outside $S$. For every extension $\mathcal{L}$ of $F$ contained in $F_S$, we write $G_S(\mathcal{L})=\Gal(F_S/\mathcal{L})$, and denote by $S(\mathcal{L})$ (resp., $S_p(\mathcal{L})$) the set of primes of $\mathcal{L}$ above $S$ (resp., $S_p$).

Let $L$ be a finite extension of $F$ contained in $F_S$. Then the fine Selmer group of $A$ over $L$ is defined by
\[ R(A/L) =\ker\left(H^1(G_S(L),\Ap)\lra \bigoplus_{v\in S(L)} H^1(L_v, \Ap)\right).\]

At first viewing, it would seem that the fine Selmer group depends on the set $S$. But we shall see that this is not so. In fact, recall that the (classical $p$-primary) Selmer group $\Sel(A/L)$ is defined by the exact sequence
\[\Sel(A/L)=\ker\left( H^1(G_S(L),\Ap)\lra \bigoplus_{v\in S(L)} H^1(L_v, A)[p^\infty]\right)\]
and it is well-known that this definition is independent of the set $S$ as long as the set $S$ contains all the primes above $p$ and the bad reduction primes of $A$ (see \cite[Chap. I, Corollary 6.6]{Mi}). Furthermore, we have a short exact sequence
\begin{equation} \label{eqn short exact} 0 \lra A(L)\ot_{\Zp}\Qp/\Zp \lra \Sel(A/L)\lra \sha(A/L)[p^\infty]\lra 0, \end{equation}
where $\sha(A/L)$ is the Tate-Shafarevich group.

 The fine Selmer group and the classical Selmer group are related by the following exact sequence.

\bl
We have an exact sequence
\[ 0\lra R(A/L) \lra \Sel(A/L) \lra \bigoplus_{v\in S_p(L)}A(L_v)\ot_{\Zp}\Qp/\Zp.\]
In particular, the definition of the fine Selmer group does not depend on the
choice of the set $S$. \el

\bpf
 See \cite[Lemma 4.1]{LMu}.
\epf

Following Wuthrich \cite{WuTS}, the fine Mordell-Weil group $\mathcal{M}(A/L)$ is defined by
\[ \mathcal{M}(A/L) = \ker\left(A(L)\ot_{\Zp}\Qp/\Zp \lra \bigoplus_{v\in S_p(L)} A(L_v)\ot_{\Zp}\Qp/\Zp \right)\]
and this fits into the following commutative diagram
\begin{equation} \label{diag def} \entrymodifiers={!! <0pt, .8ex>+} \SelectTips{eu}{}\xymatrix{
    0 \ar[r]^{} & \mM(A/L) \ar[d] \ar[r] &  A(L)\ot_{\Zp}\Qp/\Zp
    \ar[d] \ar[r] & \displaystyle\bigoplus_{v\in S_p(L)} A(L_v)\ot_{\Zp}\Qp/\Zp  \ar@{=}[d]\\
    0 \ar[r]^{} & R(A/L) \ar[r]^{} & \Sel(A/L)\ar[r] & \displaystyle\bigoplus_{v\in S_p(L)} A(L_v)\ot_{\Zp}\Qp/\Zp
     } \end{equation}
with exact rows, where the leftmost vertical map is induced by the middle vertical map. Following Wuthrich \cite{WuTS}, the fine Tate-Shafarevich group $\Zhe(A/L)$ is then defined to be
\[ \Zhe(A/L) = \coker\Big( \mM(A/L)\lra R(A/L)\Big).\]
Since the middle vertical map in diagram (\ref{diag def}) is injective, so is the leftmost vertical map. Therefore, applying the snake lemma to the diagram (\ref{diag def}), we obtain a short exact sequence
\begin{equation} \label{eqn fine short exact} 0 \lra \mM(A/L) \lra R(A/L) \lra \Zhe(A/L)\lra 0
\end{equation}
with $\Zhe(A/L)$ injecting into $\sha(A/L)$.

Let $F_\infty$ be a $\Zp$-extension of $F$, whose Galois group $\Gal(F_\infty/F)$ is denoted by $\Ga$. If $\Ga_n$ denotes the unique subgroup of $\Ga$ of index $p^n$, we write $F_n$ for the fixed field of $\Ga_n$. The fine Selmer group of $A$ over $F_\infty$ is defined to be $R(A/F_\infty) = \ilim_n R(A/F_n)$ which comes naturally equipped with a $\Zp\ps{\Ga}$-module structure. The $\Zp\ps{\Ga}$-modules $\mM(A/F_\infty)$ and $\Zhe(A/F_\infty)$ are similarly defined by taking limit of the corresponding objects over the intermediate subfields. We shall write $Y(A/F_n)$ and $Y(A/F_\infty)$ for the Pontryagin dual of $R(A/F_n)$ and $R(A/F_\infty)$ respectively. We also write $W(A/F_n)$ and $W(A/F_\infty)$ for the Pontryagin dual of $\mM(A/F_n)$ and $\mM(A/F_\infty)$ respectively. In particular, upon taking direct limit of the sequence (\ref{eqn fine short exact}) and following up by taking Pontryagin dual, we obtain
\begin{equation} \label{fine short exact} 0 \lra \Zhe(A/F_\infty)^{\vee}  \lra Y(A/F_\infty) \lra W(A/F_\infty) \lra 0.
\end{equation}

\bl \label{fg La modules}
The modules appearing in sequence $(\ref{fine short exact})$ are finitely generated over $\La$.
\el

\bpf
This is essentially well-known but for the convenience of the readers, we sketch a proof here. Since the ring $\La$ is Noetherian, it suffices to show that $Y(A/F_\infty)$ is finitely generated over $\La$, or equivalently, that $R(A/F_\infty)$ is cofinitely generated over $\La$. As $R(A/F^\cyc)$ is contained in $H^1(G_S(F_\infty),\Ap)$, we are reduced to showing that $H^1(G_S(F_\infty),\Ap)$ is cofinitely generated over $\La$. By the topological Nakayama lemma \cite[Proposition 5.3.10]{NSW}, it then suffices to show that $H^1(G_S(F_\infty),\Ap)^{\Ga}$ is cofinitely generated over $\Zp$. But since $\Ga$ has cohomological dimension one, the restriction-inflation sequence yields a surjection $H^1(G_S(F),\Ap)\tha H^1(G_S(F_\infty),\Ap)^{\Ga}$, and so we are reduced to showing that $H^1(G_S(F),\Ap)$ is cofinitely generated over $\Zp$. But the latter is a standard consequence of \cite[Proposition 8.3.20]{NSW} (for instance, see \cite[Lemma 5.5]{LMu15}).
\epf

\subsection{Control theorem for fine Selmer groups}

Retaining the settings of the previous subsection, we now state the following.

\bt \label{control theorem}
Let $A$ be an abelian variety defined over a number field $F$. Let $F_{\infty}$ be a $\Zp$-extension of $F$. Denote by $F_n$ the intermediate subfield of $F_\infty/F$ with index $|F_n:F|=p^n$. Then the restriction map
\[r_n: R(A/F_n) \lra R(A/F_{\infty})^{\Ga_n}\]
has finite kernel and cokernel which are bounded independently of $n$.
\et

Before proving the control theorem, we first establish the following lemma which is also proven in \cite[Lemma 2.0.1]{Lee}.

\bl \label{cohomology of torsion}
 Let $A$ be an abelian variety defined over $K$, where $K$ is a finite extension of either $\Q$ or $\Q_l$. Here $l$ can be any prime (possibly $=p$). Suppose that $K_\infty$ is a $\Zp$-extension of $K$ and $K_n$ is the intermediate subfield of $K_\infty$ with $|K_n:K| = p^n$. Write $G_n = \Gal(K_\infty/K_n)$. Then the group $H^1\big(G_n, A(K_\infty)[p^\infty]\big)$ is finite with order bounded independently of $n$.
\el

\bpf
We give a proof which is slightly different to that in \cite[Lemma 2.0.1]{Lee}.
Write $U : = \big(A(K_\infty)[p^\infty]\big)^\vee$. Note that this is finitely generated over $\Zp$ with a continuous action of $G:=\Gal(K_\infty/K)$. In particular, $U$ is torsion as a $\Zp\ps{G_n}$-module for every $n$. Hence we have
 \begin{equation} \label{rank} 0 = \rank_{\Zp\ps{G_n}}(U) = \rank_{\Zp}U_{G_n} - \rank_{\Zp}U^{G_n}, \end{equation}
where the second equality follows from \cite[Proposition 5.3.20]{NSW}. Now, observe  that
\[U_{G_n} = \big(A(K_n)[p^\infty]\big)^\vee, \] and this
is finite by the Mordell-Weil Theorem or Mattuck's theorem accordingly to $K$ (and hence $K_n$) being a finite extension of $\Q$ or $\Q_l$. By (\ref{rank}), so is $U^{G_n}$. In particular, we have $U^{G_n} \subseteq U[p^\infty]$. Since $U$ is finitely generated over $\Zp$, the groups $U^{G_n}$ are therefore finite and  bounded independently of $n$. But $U^{G_n}$ is precisely the Pontryagin dual of $H^1(G_n, A(K_\infty)[p^\infty])$, and so we have the conclusion of the lemma.
  \epf

We can now give the proof of Theorem \ref{control theorem}.

\bpf[Proof of Theorem \ref{control theorem}]
Consider the following commutative diagram
\[   \entrymodifiers={!! <0pt, .8ex>+} \SelectTips{eu}{}\xymatrix{
    0 \ar[r]^{} & R(A/F_n) \ar[d]^{r_n} \ar[r] &  H^{1}\big(G_{S}(F_n),A[p^\infty]\big)
    \ar[d]^{h_n} \ar[r] & \displaystyle\bigoplus_{v_n\in S(F_n)} H^1(F_{n,v_n}, \Ap) \ar[d]^{g_n=\oplus g_{n, v_n}}\\
    0 \ar[r]^{} & R(A/F_\infty)^{\Ga_n} \ar[r]^{} & H^{1}\big(G_{S}(F_\infty),A[p^\infty]\big)^{\Ga_n}\ar[r] & \displaystyle\left(\bigoplus_{w\in S(F_{\infty})} H^1(F_{\infty,w},\Ap)\right)^{\Ga_n}
     } \]
with exact rows. Since $\Ga_n$ has $p$-cohomological dimension 1, the restriction-inflation sequence tells us that $h_n$ is surjective and that $\ker h_n = H^1\big(\Ga_n, A(F_\infty)[p^\infty]\big)$. But the latter is finite with bounded order by Lemma \ref{cohomology of torsion}.

It therefore remains to show the finiteness and boundness of $\ker g_{n}$. For each $v_n$, fix a prime of $F_\infty$ above $v_n$ which is denoted by $w_n$, and write $v$ for the prime of $F$ below $v_n$. Write $\Ga_{w_n}$ for the decomposition group of $w_n$ in $\Ga$. By the Shapiro's lemma and the restriction-inflation sequence, we have
\[ \ker\Big(\bigoplus_{v_n\in S(F_n)} g_{n,v_n} \Big) = \bigoplus_{v_n\in S(F_n)} H^1\Big(\Ga_{w_n}, A(F_{\infty,v_n})[p^\infty]\Big).\]
 If $v$ is a prime of $F$ below $w_n$ such that $v$ splits completely in $F_\infty/F$, then $\Ga_{w_n}=0$ and so one has $H^1\big(\Ga_{w_n}, A(F_{\infty,w_n})[p^\infty]\big) =0$. Thus, it remains to consider the primes $v\in S$ which does not split completely in $F_\infty/F$. Since $S$ is a finite set,
the number of such possibly nonzero summands $\bigoplus H^1\big(\Ga_{w_n}, A(F_{\infty,w_n})[p^\infty]\big)$ is therefore finite and bounded independently of $n$. Hence it remains to show that each $H^1\big(\Ga_{w_n}, A(F_{\infty,w_n})[p^\infty]\big)$ is finite and bounded independently for those primes lying above $v$ which do not decompose completely in $F_\infty/F$. But this again follows from Lemma \ref{cohomology of torsion}. Thus, the proof of the theorem is completed.
 \epf

We record an immediate corollary of the control theorem.

\bc \label{control theorem cor}
Let $A$ be an abelian variety defined over a number field $F$. Let $F_{\infty}$ be a $\Zp$-extension of $F$. If $Y(A/F)$ is finite, then $Y(A/F_\infty)$ is torsion over $\La$.
\ec

It seem plausible to make the following conjecture (also see \cite{WuPhD}).

 \medskip
 \noindent
\textbf{Conjecture Y.} Let $A$ be an abelian variety defined over a number field $F$. Denote by $F_\infty$ a $\Zp$-extension with intermediate subfield $F_n$ of index $|F_n: F|=p^n$. Then $Y(A/F_\infty)$ is torsion over $\La$.

\medskip
In the event that $E$ is an elliptic curve over $\Q$ with good reduction at $p$ and $F$ is an abelian extension of $\Q$, then Conjecture Y is valid for $Y(E/F^\cyc)$ by a theorem of Kato \cite{K}. There are also some known cases for non-cyclotomic  $\Zp$-extensions (see \cite{CKL, C83, Ke, PR00, PW}).

\section{Fine Tate-Shafarevich groups over $\Zp$-extensions} \label{fine Zhe sec}

In this section, we will study the variation of the fine Tate-Shafarevich groups over a $\Zp$-extension.

\subsection{Torsionness of fine Tate-Shafarevich groups}

We begin by examining the module structure of $\Zhe(A/F_{\infty})^{\vee}$ and establishing its torsionness. We emphasis that this result does not assume the validity of Conjecture Y, and so it can be viewed as a partial evidence to Conjecture Y.

\bp \label{torsion fine TateSha}
 Let $A$ be an abelian variety defined over a number field $F$, and $F_{\infty}$ a $\Zp$-extension of $F$. Write $F_n$ for the intermediate subfield of $F_\infty/F$ with $|F_n:F| = p^n$. Suppose further that $\Zhe(A/F_n)$ is finite for every $n$. Then we have the following short exact sequence
 \[0 \lra \Zhe(A/F_{\infty})^{\vee} \lra \G\big(Y(A/F_\infty)\big) \lra \G\big(W(A/F_\infty)\big)\lra 0\]
 of
 $\La$-modules. In particular, $\Zhe(A/F_{\infty})$ is a cotorsion $\Zp\ps{\Ga}$-module.
\ep

\bpf
The second assertion follows from the short exact sequence of the proposition by Lemma \ref{G functor}(1). Therefore, it suffices to establish the said short exact sequence.
In view that $\Zhe(A/F_n)$ is finite, we may apply Lemma \ref{exact fg Zp-modules} to the dual of sequence (\ref{eqn fine short exact}) to obtain the following short exact sequence
\[0 \lra \Zhe(A/F_{n})^{\vee} \lra Y(A/F_n)[p^\infty] \lra W(A/F_n)[p^\infty]\lra 0.\]
Upon taking inverse limit, we obtain the short exact sequence
\[0 \lra \Zhe(A/F_{\infty})^{\vee} \lra \plim_n\Big(Y(A/F_n)[p^\infty]\Big) \lra \plim_n\Big(W(A/F_n)[p^\infty]\Big)\lra 0,\]
noting that $\plim_n\!^1\Zhe(A/F_{n})^{\vee}=0$ by the finiteness assumption on $\Zhe(A/F_n)$.
Hence we are now reduced to proving the isomorphisms
\[\mathfrak{G}\Big(Y(A/F_\infty)\Big) \cong \plim_n  \Big(Y(A/F_n)[p^\infty]\Big) \quad\mbox{  and }\quad \mathfrak{G}\Big(W(A/F_\infty)\Big) \cong \plim_n  \Big(W(A/F_n)[p^\infty]\Big).\]

Now, consider the exact sequence
\[ 0 \lra (\coker r_n)^{\vee} \lra Y(A/F_{\infty})_{\Ga_n}\lra Y(A/F_n)\lra  (\ker r_n)^{\vee} \lra 0.\]
 Since $(\coker r_n)^{\vee}$ is finite by Theorem \ref{control theorem}, it follows from Lemma \ref{exact fg Zp-modules} that we have an exact sequence
\[ 0 \lra (\coker r_n)^{\vee} \lra Y(A/F_{\infty})_{\Ga_n}[p^\infty]\lra Y(A/F_n)[p^\infty]\lra  (\ker r_n)^{\vee}.
\]
By Lemma \ref{limit}, this in turn yields the following exact sequence
\begin{equation} \label{res limit fine selmer seq} 0 \lra \plim_n(\coker r_n)^{\vee} \lra \plim_n\Big(Y(A/F_{\infty})_{\Ga_n}[p^\infty]\Big)\lra \plim_nY(A/F_n)[p^\infty]\lra  \plim_n(\ker r_n)^{\vee}.
\end{equation}
 On the other hand, it follows from the definition of the restriction map
 \[r_n: R(A/F_n) \lra R(A/F_{\infty})^{\Ga_n} \]
 that one has $\ilim_n \ker r_n = \ilim_n\coker r_n = 0$ which in turn implies that
$\plim_n (\ker r_n)^{\vee} = \plim_n(\coker r_n)^{\vee}= 0$. Combining these observations with the exact sequence (\ref{res limit fine selmer seq}), we obtain
\[ \mathfrak{G}\Big(Y(A/F_\infty)\Big) \cong \plim_n  \Big(Y(A/F_n)[p^\infty]\Big), \]
 as required.
Now consider the following commutative diagram
\[   \entrymodifiers={!! <0pt, .8ex>+} \SelectTips{eu}{}\xymatrix{
    0 \ar[r]^{} & \mM(A/F_n) \ar[d]^{t_n} \ar[r] &  R(A/F_n)
    \ar[d]^{r_n} \ar[r] & \Zhe(A/F_n)  \ar[d]^{z_n} \ar[r] & 0\\
    0 \ar[r]^{} & \mM(A/F_\infty)^{\Ga_n} \ar[r]^{} & R(A/F_\infty)^{\Ga_n}\ar[r] & \Zhe(A/F_\infty)^{\Ga_n} &
     } \]
with exact rows. Since $\Zhe(A/F_n)$ is assumed to be finite for all $n$, so is $\ker z_n$. Combining this with Theorem \ref{control theorem}, we see that
$\ker t_n$ and $\coker t_n$ are finite for all $n$. We can now proceed similarly as above to conclude that
\[\mathfrak{G}\Big(W(A/F_\infty)\Big) \cong \plim_n  \Big(W(A/F_n)[p^\infty]\Big). \]
This completes the proof of the proposition.
\epf

\bc \label{torsion fine TateSha2}
 Retain the setting of Proposition \ref{torsion fine TateSha}. Suppose that $\Zhe(A/F_n)$ is finite for every $n$ and that $\big(A(F_\infty)\ot_{\Zp}\Qp/\Zp\big)^{\vee}$ is finitely generated over $\Zp$. Then we have $\G\big(W(A/F_\infty)\big)=0$ and an isomorphism
 \[\Zhe(A/F_{\infty})^{\vee} \cong\G\big(Y(A/F_\infty)\big)\]
 of
 torsion $\La$-modules.
\ec

\bpf
In view of the hypothesis, it follows from the result of Lee \cite[Theorem 2.1.2]{Lee} that we have an injection
\[ \big(A(F_\infty)\ot_{\Zp}\Qp/\Zp\big)^{\vee} \hookrightarrow \bigoplus_{j=1}^s \La/\xi_{n_j} \]
with finite cokernel, where $\xi_{n_j}$ is certain $p^{n_j}$th cyclotomic polynomials. Since $\big(A(F_\infty)\ot_{\Zp}\Qp/\Zp\big)^{\vee}$ surjects onto $W(A/F_\infty)$, we see that $W(A/F_\infty)$ is also pseudo-isomorphic to a direct summand of modules of the form $\La/\xi_{n_j}$. It then follows from Proposition \ref{G functor}(3) that $\G\big(W(A/F_\infty)\big)=0$. Putting this latter observation into the short exact sequence of Proposition \ref{torsion fine TateSha}, we obtain the conclusion of the corollary.
\epf

\br
In proving $\G\Big(\big(A(F_\infty)\ot_{\Zp}\Qp/\Zp)\big)^{\vee}\Big)=0$ in \cite[Theorem 2.1.2]{Lee}, a crucial observation used by Lee is that $A(F_n)\ot_{\Zp}\Qp/\Zp$ is $p$-divisible. Unfortunately, this property may not hold for the fine Mordell-Weil group $\mM(A/F_n)$, and so we cannot apply the approach there to prove $\G\big(W(A/F_\infty)\big)=0$ directly and unconditionally. \er

\subsection{Growth of fine Tate-Shafarevich groups}

We give the following growth formula for the $p$-torsion subgroup of $Y(A/F_n)$ which is an immediate consequence of our control theorem. Note that this result does not require the validity of Conjecture Y,

\bp \label{growthfineselmer}
 Let $A$ be an abelian variety defined over a number field $F$. Let $F_{\infty}$ be a $\Zp$-extension of $F$ and $F_n$ the intermediate subfield of $F_\infty$ with $|F_n:F| = p^n$. Then we have
 \[ e\big(Y(A/F_n)\big) = \mu\Big(\G\big(Y(A/F_\infty)\big)\Big)p^n +\la\Big(\G\big(Y(A/F_\infty)\big)\Big)n + O(1).\]
\ep

\bpf
 This is an immediate consequence of Lemma \ref{tech estimate} and Theorem \ref{control theorem}.
\epf

As mentioned in the introduction, the fine Mordell-Weil group $\mM(A/F_n)$ is not necessarily $p$-divisible, and so we \textit{may not} have an equality $\Zhe(A/F_n) = Y(A/F_n)[p^\infty]$ in general. Therefore, the growth formula in Proposition \ref{growthfineselmer} does not transfer directly to the fine Tate-Shafarevich groups.

We do at the very least have the following lower bound and upper bound for  the growth of the fine Tate-Shafarevich groups.

\bp \label{growth bound}
 Let $A$ be an abelian variety defined over a number field $F$. Let $F_{\infty}$ be a $\Zp$-extension of $F$ and $F_n$ the intermediate subfield of $F_\infty$ with $|F_n:F| = p^n$. Suppose further that  $\Zhe(A/F_n)$ is finite for every $n$.
    Then we have
 \[ \mu\Big(\Zhe(E/F_\infty)^\vee\Big)p^n+  \la\Big(\Zhe(E/F_\infty)^\vee\Big)n + O(1) \leq e\big(\Zhe(A/F_n)\big) \hspace{2in}\]\[\hspace{2in} \leq \mu\Big(\G\big(Y(A/F_\infty)\big)\Big)p^n +\la\Big(\G\big(Y(A/F_\infty)\big)\Big)n + O(1).\]
\ep

\bpf
As seen in the proof of Proposition \ref{torsion fine TateSha}, one has a short exact sequence
 \[0 \lra \Zhe(A/F_{n})^{\vee} \lra Y(A/F_n)[p^\infty] \lra W(A/F_n)[p^\infty]\lra 0\]
 which in turn implies that
 \[  e\big(\Zhe(A/F_n)\big)= e\big(\Zhe(A/F_n)^{\vee}\big) = e\big(Y(A/F_n)\big) - e\big(W(A/F_n)\big) .\]
  By Proposition \ref{growthfineselmer}, we have
 \[ e\big(Y(A/F_n)\big) = \mu\Big(\G\big(Y(A/F_\infty)\big)\Big)p^n +\la\Big(\G\big(Y(A/F_\infty)\big)\Big)n + O(1). \]
 On the other hand, it follows from Lemma \ref{tech estimate} that
 \[ e\big(W(A/F_n)\big) \leq \mu\Big(\G\big(W(A/F_\infty)\big)\Big)p^n + \la\Big(\G\big(W(A/F_\infty)\big)\Big)n + O(1). \]
 Combining these estimates and taking Proposition \ref{torsion fine TateSha} into account, we obtain the proposition.
  \epf

It is naturally to ask whether one can derive a precise asymptotic formula for the growth of the fine Tate-Shafarevich groups under appropriate extra assumptions.
The following two results (Theorems \ref{growth} and \ref{growth2}) provide sufficient conditions for this.

\bt \label{growth}
 Let $A$ be an abelian variety defined over a number field $F$. Let $F_{\infty}$ be a $\Zp$-extension of $F$ and $F_n$ the intermediate subfield of $F_\infty/F$ with $|F_n:F| = p^n$. Suppose that $A(F_\infty)$ is a finitely generated abelian group and that $\Zhe(A/F_n)$ is finite for each $n$.
      Then we have
 \[e\big(\Zhe(A/F_n)\big) = \mu\Big(\Zhe(A/F_\infty)^\vee\Big)p^n + \la\Big(\Zhe(A/F_\infty)^\vee\Big)n + \nu\]
 for $n\gg0$.
\et

\bpf
As seen in the proof of Proposition \ref{growth bound}, we have
\[  e\big(\Zhe(A/F_n)\big)= e\big(\Zhe(A/F_n)^{\vee}\big) = e\big(Y(A/F_n)\big) - e\big(W(A/F_n)\big).\]
  By Proposition \ref{growthfineselmer} and Corollary \ref{torsion fine TateSha2}, we have
 \[ e\big(Y(A/F_n)\big) = \mu\Big(\Zhe(A/F_\infty)^{\vee}\Big)p^n + \la\Big(\Zhe(A/F_\infty)^{\vee}\Big)n + c \]
 for $n\gg 0$. It therefore remains to show that $e\big(W(A/F_n)\big)$ is bounded independently of $n$. As seen in the proof of Proposition \ref{torsion fine TateSha}, the module $W(A/F_\infty) = \plim W(A/F_n)$ satisfies the hypothesis of Lemma \ref{tech estimate}. Therefore, we may combine Corollary \ref{torsion fine TateSha2} with Lemma \ref{tech estimate} to obtain the boundness of $e\big(W(A/F_n)\big)$ and this completes the proof of the theorem.
\epf

\bt \label{growth2}
 Let $A$ be an abelian variety defined over a number field $F$ with potentially good ordinary reduction at all primes of $F$ above $p$. Let $F_{\infty}$ be a $\Zp$-extension of $F$ and $F_n$ the intermediate subfield of $F_\infty/F$ with $|F_n:F| = p^n$. Suppose that $\sha(A/F_n)[p^\infty]$ is finite for every $n$.
    Then we have
 \[e\big(\Zhe(A/F_n)\big) = \mu\Big(\Zhe(A/F_\infty)^\vee\Big)p^n + \la\Big(\Zhe(A/F_\infty)^\vee\Big)n + \nu\]
 for $n\gg0$.
\et

\bpf
As seen in the proof of Theorem \ref{growth}, we have
\[  e\big(\Zhe(A/F_n)\big)= e\big(\Zhe(A/F_n)^{\vee}\big) = e\big(Y(A/F_n)\big) - e\big(W(A/F_n)\big).\]
  By virtue of Propositions \ref{torsion fine TateSha} and \ref{growthfineselmer}, we are reduced to showing that
  \[ e\big(W(A/F_n)\big) = \mu\Big(\G\big(W(A/F_\infty)\big)\Big)p^n + \la\Big(\G\big(W(A/F_\infty)\big)\Big)n + d \]
   for $n\gg 0$. Now taking Lemma \ref{tech estimate} into account, it suffices to show that the kernel and cokernel of the map
   \[ t_n: \mM(A/F_n)\lra \mM(A/F_\infty)^{\Ga_n}\]
   are finite and bounded independently of $n$. Since $\sha(A/F_n)$ is assumed to be finite for all $n$, so is $\Zhe(A/F_n)$. Therefore, as seen in the proof of Proposition \ref{torsion fine TateSha}, we see that $\ker t_n$ and $\coker t_n$ is finite. It therefore remains to show that they are bounded independently on $n$. Plainly, the kernel of $t_n$ is contained in the kernel of
   \[ A(F_n)\ot\Qp/\Zp \lra \Big(A(F_\infty)\ot\Qp/\Zp\Big)^{\Ga_n}\]
   which in turn is contained in the kernel of
   \[ H^1(G_S(F_n),\Ap) \lra H^1\big(G_S(F_\infty),\Ap\big)^{\Ga_n}.\]
   But this is finite and bounded independently of $n$ by Lemma \ref{cohomology of torsion}.

   We now show that $t_n$ has finite cokernel which is bounded independently of $n$. Consider the following commutative diagram
   \[   \entrymodifiers={!! <0pt, .8ex>+} \SelectTips{eu}{}\xymatrix{
    0 \ar[r]^{} & A(F_n)\ot_{\Zp}\Qp/\Zp\ar[d]^{a_n} \ar[r] &  \Sel(A/F_n)
    \ar[d]^{s_n} \ar[r] & \sha(A/F_n)[p^\infty]  \ar[d] \ar[r] & 0\\
    0 \ar[r]^{} & \big(A(F_\infty)\ot_{\Zp}\Qp/\Zp\big)^{\Ga_n} \ar[r]^{} & \Sel(A/F_\infty)^{\Ga_n}\ar[r] & \big(\sha(A/F_\infty)[p^\infty]\big)^{\Ga_n} &
     } \]
with exact rows. By hypothesis, the control theorem for the classical Selmer group tells us that $\coker s_n$ is finite for each $n$ (cf. \cite[Proposition 5.1]{G03}). Since $\sha(A/F_n)[p^\infty]$ is finite by hypothesis, it follows that $\coker a_n$ is finite for each $n$. As observed by Lee \cite[Lemma 4.2.4]{Lee}, the finiteness of $\coker a_n$ in turn implies that they are automatically bounded independently of $n$.
Now consider the following commutative diagram
\[   \entrymodifiers={!! <0pt, .8ex>+} \SelectTips{eu}{}\xymatrix{
    0 \ar[r]^{} & \mM(A/F_n) \ar[d]^{t_n} \ar[r] &  A(F_n)\ot_{\Zp}\Qp/\Zp\ar[d]^{a_n} \ar[r] & \displaystyle\bigoplus_{v_n\in S_p(F_n)}A(F_{n, v_n})\ot_{\Zp}\Qp/\Zp \ar[d]^{b_n} \\
    0 \ar[r]^{} & \mM(A/F_\infty)^{\Ga_n} \ar[r]^{} & \big(A(F_\infty)\ot_{\Zp}\Qp/\Zp\big)^{\Ga_n}\ar[r] & \displaystyle\Big(\bigoplus_{w\in S_p(F_\infty)}A(F_{\infty,w})\ot_{\Zp}\Qp/\Zp\Big)^{\Ga_n}
     } \]
with exact rows. In view of the boundedness of $\coker a_n$, it remains to show that $\ker b_n$ is finite and bounded independently of $n$. But $\ker b_n$ is contained in the kernel of the following map
\[ \displaystyle\bigoplus_{v_n\in S_p(F_n)}H^1(F_{n,v_n},\Ap) \lra \left(\bigoplus_{w\in S_p(F_\infty)}H^1(F_{\infty,w},\Ap)\right)^{\Ga_n}, \]
and we have seen that the latter is finite and bounded in the proof of Theorem \ref{control theorem}. This completes the proof of the theorem.
\epf

We end the section with the natural question.

\medskip \noindent
\textbf{Question 1.} Does one always have  \[e\big(\Zhe(A/F_n)\big) = \mu\Big(\Zhe(A/F_\infty)^\vee\Big)p^n + \la\Big(\Zhe(A/F_\infty)^\vee\Big)n + \nu\]
 for $n\gg0$?

\medskip
In fact, we shall see in the next section that we do not have an answer even restricting to a specific class of $\Zp$-extensions.

\section{Further comments} \label{further comments}

\subsection{Cyclotomic $\Zp$-extension}
Write $F^\cyc$ for the cyclotomic $\Zp$-extension of $F$.
It has been conjectured that $A(F^\cyc)$ is a finitely generated abelian group (see \cite{G99, Maz}; also see \cite{Lee, LP}). Therefore, assuming this conjecture, one expects to have an asymptotic formula for the fine Shafarevich groups by Theorem \ref{growth2}(i). We record the following important case, where we do have such formula unconditionally.

\bp
Let $E$ be an elliptic curve defined over $\Q$ with good reduction at $p$, and $F$ an abelian extension of $\Q$. Denote by $F^\cyc$ the cyclotomic $\Zp$-extension with intermediate subfield $F_n$ of $|F_n: F|=p^n$. Suppose that $\Zhe(E/F_n)$ is finite for all $n$. Then we have
 \[e\big(\Zhe(E/F_n)\big) = \mu\Big(\Zhe(E/F^\cyc)^\vee\Big)p^n + \la\Big(\Zhe(E/F^\cyc)^\vee\Big)n + \nu\]
 for $n\gg0$.
\ep

\bpf
Under the hypothesis of the proposition, a well-known result of Kato \cite[Theorem 14.4]{K} and Rohrlich \cite{Ro} asserts that $E(F^\cyc)$ is finitely generated. The conclusion is then immediate from Theorem \ref{growth}.
\epf

In a recent paper of Lei-Ponisnet \cite{LP}, they have given a sufficient condition for the finite generation of $A(F^\cyc)$ for an abelian variety with good supersingular reduction at all primes above $p$. Therefore, one can also obtain a nice growth formula for the fine Tate-Shafarevich groups in their context.
We should mention that the conclusion of the corollary was also stated in \cite[Section 8]{WuTS} (but without a proof). Coates and Sujatha have further conjectured that $Y(A/F^\cyc)$ should be finitely generated over $\Zp$ (see \cite[Conjecture A]{CS05}). In view of their conjecture, we expect $\Zhe(A/F^\cyc)^\vee$ to have trivial $\mu$-invariant. Therefore, it seems plausible to make the following conjecture.

 \medskip
 \noindent
\textbf{Conjecture Z.} Let $A$ be an abelian variety defined over a number field $F$. Denote by $F^\cyc$ the cyclotomic $\Zp$-extension with intermediate subfield $F_n$ of index $|F_n: F|=p^n$. Suppose that $\Zhe(A/F_n)$ is finite for all $n$. Then one has
 \[e\big(\Zhe(A/F_n)\big) = \la\Big(\Zhe(A/F^\cyc)^\vee\Big)n + \nu\]
 for $n\gg0$.

 \medskip
The above conjecture was also stated in \cite[Conjecture 8.2]{WuTS} for an elliptic curve (also see \cite{WuPhD}).

\br
We mention some evidences on the conjecture.
\begin{enumerate}
\item[$(i)$] Suppose that $E$ is an elliptic curve defined over $\Q$ and $F$ is an abelian extension of $\Q$ with $E(F)[p]\neq 0$, then it has been shown that $Y(E/F^\cyc)$ is finitely generated over $\Zp$ (cf. \cite[Corollary 3.6]{CS05}). Therefore, the above conjecture is valid.

\item[$(ii)$] Suppose that $A$ is an abelian variety with good ordinary reduction at all primes above $p$ and that $\Sel(A/F^\cyc)^{\vee}$ is a torsion $\La$-module with trivial $\mu$-invariant. Since $(A(F^\cyc)\ot_{\Zp}\Qp/\Zp)^{\vee}$ and $Y(A/F^\cyc)$ are quotients of $\Sel(A/F^\cyc)^{\vee}$, we have that $A(F^\cyc)$ is finitely generated and $Y(A/F^\cyc)$ has trivial $\mu$-invariant. Therefore, the conjecture holds in this case.

Furthermore, if $L$ is a finite $p$-extension of $F$, then $\Sel(A/L^\cyc)^{\vee}$ is also a torsion $\La$-module with trivial $\mu$-invariant (cf. \cite{HM}), and so we also have the validity of the conjecture for these extensions.
Therefore, this can provide many examples of where one has the validity of the conjecture. For instance, for $p=5$, and $E$ is the elliptic curve $y^2 + y = x^3 -x^2$. It is well-known that $\Sel(E/\Q(\mu_5)^\cyc)^\vee=0$ (see \cite[Theorem 5.4]{CS10}; we thank the anonymous referee for reminding us of this fact). In particular, it is a torsion $\La$-module with trivial $\mu$-invariant. Hence the same can be said for $\Sel(E/L^\cyc)^\vee$ for any finite $p$-extension $L$ of $\Q(\mu_5)$. Examples of such $L$ are $\Q(\mu_{5^i}, \sqrt[5^{j_1}]{m_1},...,\sqrt[5^{j_r}]{m_r})$ for $i\geq j_k$ and $m_k$ a $p$-powerfree integer.

\item[$(iii)$] Let $E$ be an elliptic curve defined over $\Q$ with good supersingular reduction at $p$. The Pontryagin dual of the plus-minus Selmer groups $\Sel^{\pm}(E/\Q^{\cyc})$ in the sense of Kobayashi \cite{Kob} (we refer readers therein for the precise definitions of these groups) is $\La$-torsion (see \cite[Theorem 1.2]{Kob}). Furthermore, it is conjectured that these signed Selmer groups have trivial $\mu$-invariants (see discussion before \cite[Corollary 10.10]{Kob} and \cite{Po}) and this have been numerical verified for many elliptic curves \cite{PoCom}. Therefore, these examples gives many cases where one has the growth formula for the fine Tate-Shafarevich groups as conjectured.
\end{enumerate}
\er

Over the cyclotomic $\Zp$-extension of $\Q$, Wuthrich even questioned whether the fine Tate-Shafarevich groups can have bounded growth (see \cite[Question 8.3]{WuTS}; also see \cite{C12, WuPhD}). We shall have little to say on this. But
this naturally leads us to the next question (we thank Christian Wuthrich for pointing this out) which we have no answer at this point of writing.

\medskip \noindent
\textbf{Question 2.} Does there exist an example of $\Zhe(A/F^\cyc)$ being infinite?
\medskip

\subsection{Elliptic curves over certain $\Zp$-extensions of an imaginary quadratic field}

In this subsection, we consider certain non-cyclotomic $\Zp$-extensions. Throughout the discussion here, $p$ is taken to be a prime $\geq 5$.
Let $E$ be an elliptic curve defined over $\Q$ with good reduction at the prime $p$ and $F$ an imaginary quadratic field of $\Q$. We first recall the following conjecture of Mazur \cite{Maz83}.

 \medskip \noindent
\textbf{Growth Number
Conjecture.}~(Mazur) The Mordell-Weil rank of $E$ stays bounded along
any $\Zp$-extension of $F$, unless the extension is anticyclotomic and the root
number is negative.

\medskip
This conjecture was stated by Mazur in \cite[Growth Number
Conjecture]{Maz83} for an ellptic curve with good ordinary reduction at $p$. A variant of this conjecture was stated and studied in \cite[Conjecture 1.2]{LSp} for an elliptic curve with good supersingular reduction at $p$. In view of this conjecture, one would expect to have a precise growth formula for the fine Tate-Shafarevich groups over infinitely many of such $\Zp$-extensions. We are therefore led to the following question, which is a special case of Question 1.

\medskip \noindent
\textbf{Question 3.} Denote by $F^\ac$ the anticyclotomic $\Zp$-extension of the imaginary quadratic field $F$ and write $F_n$ for the intermediate subfield of $F^\ac$ with $|F_n: F|=p^n$.  Does one always have  \[e\big(\Zhe(A/F_n)\big) = \mu\Big(\Zhe(A/F^\ac)^\vee\Big)p^n + \la\Big(\Zhe(A/F^\ac)^\vee\Big)n + \nu\]
 for $n\gg0$?
\medskip

Now if the elliptic curve $E$ has good ordinary reduction at $p$, then Theorem \ref{growth2} applies (modulo finiteness of $\sha(E/F_n)$). Therefore, what remains is the case of an elliptic curve with supersingular reduction at $p$ and when the root number is negative. We do not have an answer to this and hope to revisit this in a subsequent work.

For the remainder of the subsection, we mention some cases, where we have a precise growth formula and raise some question on these cases. Clearly, if the elliptic curve has good ordinary reduction at $p$, Theorem \ref{growth2} applies to yield the formula (assuming finiteness of $\sha(E/F_n)$). Therefore, our presentation below will focus more from the point of view of Theorem \ref{growth}. As a start, we have the following observation.

\bp \label{ac}
Let $E$ be an elliptic curve of squarefree conductor $N$ which has good reduction at the prime $p\geq 5$. Let $F$ be an imaginary quadratic field with discriminant coprime to $pN$. Write $N = N^+N^-$ with $N^+$ divisible
only by primes which are split in $F/\Q$ and $N^-$ divisible only by inert primes. Suppose that $N^-$ has an odd number of prime divisor. Let $F^\ac$ be the anticyclotomic $\Zp$-extension. Write $F_n$ for the intermediate subfield of $F^\ac$ with $|F_n: F|=p^n$. Suppose that $\Zhe(E/F_n)$ is finite for all $n$. In the event that $E$ has supersingular reduction at $p$, assume further that $p$ is split in $F/\Q$ and that each prime of $F$ above $p$ is totally ramified in $F^\ac/F$. Then we have
 \[e\big(\Zhe(E/F_n)\big) =  \la\Big(\Zhe(E/F^\ac)^\vee\Big)n + \nu\]
 for $n\gg0$.
\ep

\bpf
 By the result of Pollack-Weston \cite[Theorems 1.1 and 1.3]{PW} (also see \cite{Vat}), the dual Selmer group (resp., dual signed Selmer groups) of $E$ over $F^\ac$ is torsion over $\La$ with trivial $\mu$-invariants, when $E$ has good ordinary reduction at $p$ (resp, good supersingular reduction at $p$). This yields the finite generation of $E(F^\ac)$ (also see \cite{BD, DI, Vat}). Finally, as the fine Selmer group sits in every Selmer groups, we also have the vanishing of the $\mu$-invariant of $Y(E/F^\ac)$. Putting these information into Theorem \ref{growth}, we have the asymptotic formula as asserted.
\epf

We now consider a complementary situation of the preceding proposition. Suppose that $E$ is an elliptic curve of squarefree conductor $N$ which has good reduction at the prime $p\geq 5$. Let $F$ be an imaginary quadratic field with discriminant coprime to $pN$ and such that all the prime divisor of $pN$ splits completely in $F/\Q$. Matar \cite[Conjecture B]{Mat18} has conjectured that $Y(E/F^\ac)$ is finitely generated over $\Zp$ and provided some evidence to his conjecture \cite[Theorem 4.1]{Mat18}.

In view of Proposition \ref{ac} and this conjecture of Matar, one might ask the following question.

\medskip \noindent
\textbf{Question 4.} Does one always have $\mu\Big(\Zhe(E/F^\ac)^\vee\Big) =0$?
\medskip

We do not have an answer to our question at this point of writing. In fact, we do not even have an answer to a modest variant of the above question. Before stating this, we record the following.

\bp \label{ac cm}
Let $E$ be an elliptic curve with complex multiplication given by the ring of integers of an imaginary quadratic field $F$. Suppose that $p$ is a prime $\geq 5$ which split completely in $F/\Q$. Write $F_n$ for the intermediate subfield of $F^\ac/F$ with $|F_n: F|=p^n$. Suppose that either of the following statements is valid.
\begin{enumerate}
\item[$(a)$] $\Zhe(E/F_n)$ is finite for all $n$ and that the root number of $E/\Q$ is $+1$.
\item[$(b)$]  $\sha(E/F_n)$ is finite for all $n$.
\end{enumerate}
Then we have
 \[e\big(\Zhe(E/F_n)\big) = \mu\Big(\Zhe(E/F^\ac)^\vee\Big)p^n + \mu\Big(\Zhe(E/F^\ac)^\vee\Big)n + \nu\]
 for $n\gg0$.
\ep

\bpf
Since the prime $p$ splits completely in $F/\Q$, $E$ has good ordinary reduction at all primes above $p$. So the conclusion follows from Theorem \ref{growth2} under assumption (b). To see that one has the asserted growth formula under (a), recall that under the root number hypothesis, a well-known theorem of Greenberg \cite[Theorem 3]{G83} asserts that $E(F^\ac)$ is finitely generated, and so one can apply Theorem \ref{growth}.
\epf

We can now ask the following modest variant of Question 4.

\medskip \noindent
\textbf{Question 5.} Retain the setting of Proposition \ref{ac cm}. Does one have $\mu\Big(\Zhe(E/F^\ac)^\vee\Big) =0$?
\medskip

In the case of root number +1, the work of Finis \cite{Fin} shows that the $\mu$-invariant of the anticyclotomic $p$-adic $L$-function can possibly be nonzero (We thank Ming-Lun Hsieh for explaining this to us). In view of the anticyclotomic main conjecture (for instance, see \cite{AH}), the Pontryagin dual of the classical Selmer group may possibly have nonzero $\mu$-invariant. However, since the fine Selmer group is smaller than the Selmer group, this does not necessarily imply that the fine Selmer group has non-trivial $\mu$-invariant. For the case of the root number being $-1$, we have even less information, since the classical Selmer group is not cotorsion over $\La$ anymore (see \cite{AH}), and so it does not seem easy to extract information on the structure of the dual fine Selmer group from the non-torsion classical Selmer groups.

We end by mentioning another class of $\Zp$-extensions of $F$ coming from CM elliptic curves.

\bp \label{cm}
Let $E$ be an elliptic curve with complex multiplication given by the ring of integers of an imaginary quadratic field $F$. Suppose that $p$ is a prime $\geq 5$ which split completely in $F/\Q$, say $p = \mathfrak{p}\bar{\mathfrak{p}}$, where $\mathfrak{p}$ is a prime of $F$ above $p$. Let $F_{\mathfrak{p}^\infty}$ be the unique $\Zp$-extension of $F$ unramified outside $\mathfrak{p}$. Write $F_n$ for the intermediate subfield of $F_{\mathfrak{p}^\infty}/F$ with $|F_n: F|=p^n$. Suppose that $\Zhe(E/F_n)$ is finite for all $n$. Then we have
 \[e\big(\Zhe(E/F_n)\big) =  \la\Big(\Zhe(E/F_{\mathfrak{p}^\infty})^\vee\Big)n + \nu\]
 for $n\gg0$.
\ep

\bpf
By a classical result of Coates \cite[Theorem 16]{C83} (or see \cite[Chap IV, Corollary 1.8]{deS}), $E(F_{\mathfrak{p}^\infty})$ is a finitely generated abelian group. Furthermore, via the results of Gillard \cite{Gi} and Schneps \cite{Sch}, one can show that the dual strict Selmer group of $E$ over $F_{\mathfrak{p}^\infty}$ is torsion over $\La$ with trivial $\mu$-invariant, which in turn implies that the $\mu$-invariant of $Y(E/F_{\mathfrak{p}^\infty})$ is trivial (see \cite[Proposition 5.1]{LSu2} for details). Hence we may apply Theorem \ref{growth} to obtain the conclusion of the proposition.
\epf

\end{document}